\documentclass[12pt]{article}

\usepackage{anysize}
\marginsize{2cm}{2cm}{1cm}{2cm}
\usepackage{amsmath}
\usepackage{amssymb}

\newcommand{\proof}{\par\noindent{\it Proof.\ \ }}
\newcommand{\qed}{\ifmmode\square\else\nolinebreak\hfill
$\Box$\fi\par\vskip12pt}

\newcommand\Aut{{\sf Aut}} 
\newcommand\Ker{{\sf Ker}} 
\newcommand\Sym{{\sf Sym}}
\newcommand\GF{{\sf GF}}
\newcommand\PGammaL{{\sf P\Gamma L}} 
\newcommand\PSL{{\sf PSL}}

\newcommand\AGL{{\sf AGL}}

\newcommand{\norml}{\vartriangleleft}
\newcommand\la{\langle}
\newcommand\ra{\rangle}
\newcommand\calP{{\mathcal P}}
\newcommand\Ga{\Gamma}
\newcommand\sig{\sigma}
\newcommand\diam{{\sf diam}}
\newcommand\girth{{\sf girth}}
\newcommand\C{{\bf C}} 
\newcommand\Cos{{\sf Cos}}

\newcommand\N{{\bf N}}
\newcommand\Z{{\bf Z}}
 
\newcommand\ga{g_\alpha}
\newcommand\ua{u_\alpha}
\newcommand\dist{{\sf d}}
\DeclareMathOperator{\Wr}{wr}

\newtheorem{theorem}{Theorem}[section]%
\newtheorem{lemma}[theorem]{Lemma}%
\newtheorem{corollary}[theorem]{Corollary}%
\newtheorem{proposition}[theorem]{Proposition}%
\newtheorem{definition}[theorem]{Definition}%
\newtheorem{construction}[theorem]{Construction}%
\newtheorem{question}[theorem]{Question}%
\newtheorem{remark}[theorem]{Remark}%

\title
{An infinite family of biquasiprimitive $2$-arc transitive cubic graphs}
\author{Alice Devillers,
 Michael Giudici, 
Cai Heng Li and Cheryl E.~Praeger\footnote{The paper forms part of Australian Research Council Federation Fellowship FF0776186 held by the fourth author. The first author is supported by UWA as part of the Federation Fellowship project and the second author is supported by an Australian Research Fellowship. }\footnote{emails: alice.devillers@uwa.edu.au, michael.giudici@uwa.edu.au, cai.heng.li@uwa.edu.au, cheryl.praeger@uwa.edu.au}\\
Centre for the Mathematics of Symmetry and Computation\\
School of Mathematics and Statistics\\
The University of Western Australia\\
35 Stirling Highway\\
Crawley WA 6009\\
Australia\\
 }
\date\today

\begin{document}

\maketitle
 \begin{abstract}
A new infinite family of bipartite cubic 3-arc transitive graphs is constructed and studied. They provide the first known examples admitting a 2-arc transitive vertex-biquasiprimitive group of automorphisms for which the index two subgroup fixing each half of the bipartition is not quasiprimitive on either bipartite half.
\end{abstract}
\section{Introduction}

The study of cubic $s$-arc-transitive graphs goes back to the seminal papers of Tutte \cite{Tutte1,Tutte2} who showed that $s\leq 5$. More generally,  Weiss \cite{8-arc-trans} proved that $s\leq 7$ for graphs of larger valency. In \cite{quasi}, the last author introduced a global approach to the study of $s$-arc-transitive graphs. 

Given a connected graph $\Gamma$ with an $s$-arc-transitive group $G$ of automorphisms, if $G$ has a nontrivial normal subgroup $N$ with at least three orbits on vertices, then $G$ induces an unfaithful but $s$-arc-transitive action on the normal quotient $\Gamma_N$ (defined in Section \ref{sec:prelim}). The important graphs to study are then those with no ``useful'' normal quotients, that is, those for which all nontrivial normal subgroups of $G$ have at most two orbits on vertices. A transitive permutation group for which all nontrivial normal subgroups are transitive is called \emph{quasiprimitive}, while if the group is not quasiprimitive and all nontrivial normal subgroups have at most two orbits we call it \emph{biquasiprimitive}. Thus the basic graphs to study are those which are $(G,s)$-arc transitive and $G$ is either quasiprimitive or biquasiprimitive on vertices.

Now suppose that our graph $\Gamma$ were bipartite. Then the \emph{even subgroup} $G^+$ (the subgroup generated by the vertex stabilisers $G_v$ for all $v\in V\Gamma$) has index 2 in $G$ and is transitive on each of the two bipartite halves of $\Gamma$ (see, for example, \cite[Proposition 1]{DjoMil}). Since $G^+$ is vertex-intransitive, $G$ is not vertex-quasiprimitive and so the basic bipartite graphs are those where $G$ is biquasiprimitive on vertices. The actions of such groups were investigated in \cite{bqp,bip}. However, when $G$ is biquasiprimitive it may still be possible to find a meaningful quotient of the graph. The subgroup $G^+$ is what is called  locally transitive on $s$-arcs (see Section \ref{sec:prelim} for precise definition and \cite{GLP} for an analysis of such graphs).  If $G^+$ is not quasiprimitive on each bipartite half (note the two actions of $G^+$ are equivalent) then we can form a $G^+$-normal quotient and obtain a new (smaller) locally $s$-arc-transitive graph. The existence of a 2-arc transitive graph with such a group has been regarded as `problematic' (see \cite[Section 4]{bqp}). The main result of this paper is that there do indeed exist $(G,2)$-arc transitive graphs such that $G$ is biquasiprimitive but $G^+$ is not quasiprimitive on each bipartite half.

\begin{theorem}\label{thm:broadbrush}
There exist infinitely many connected bipartite $(G,2)$-arc transitive graphs $\Gamma$ of valency $3$, where $G\leq \Aut(\Gamma)$, such that $G$ is biquasiprimitive on vertices but $G^+$ is not quasiprimitive on either bipartite half.
\end{theorem}

Such permutation groups $G$ were described in detail in \cite[Theorem 1.1(c)(i)]{bqp} (see Corollary \ref{cor:(c)(i)}) and this theorem gives the first examples of 2-arc-transitive graphs admitting such an automorphism group. (Our graphs are actually $3$-arc transitive, but only $(G,2)$-arc-transitive.)
We also  provide an infinite family of $(G,1)$-arc-transitive graphs where $G$ is biquasiprimitive on vertices but $G^+$ is not quasiprimitive on each orbit (Construction \ref{s=1}). The full automorphism group $A$ of these graph is 2-arc-transitive but $A^+$ is quasiprimitive on each bipartite half.

Graphs which are $s$-arc transitive are also $s$-distance transitive, provided their diameter is at least $s$. Such graphs were studied in \cite{LDT} where $(G,s)$-distance transitive bipartite graphs with $G$  biquasiprimitive on vertices but $G^+$ not quasiprimitive on each bipartite half were referred to as $G$-basic but not $G^+$-basic (see \cite[Proposition 6.3]{LDT}). Our infinite family of graphs shows that connected 2-distance transitive graphs with such an automorphism group do indeed exist and so this answers Question 6.4 of \cite{LDT} in the affirmative for $s=2$.

We prove Theorem \ref{thm:broadbrush} by constructing and analysing a new infinite family of finite bipartite $(G,2)$-arc transitive graphs $\Gamma(f,\alpha)$ of valency 3, where $f$ is a positive integer and $\alpha$ lies in the Galois field $\GF(2^f)$, see Construction~\ref{cons}.  The group $G\leqslant \Aut(\Gamma(f,\alpha))$ depends only on $f$, has order $2^{2f+1}(2^{2f}-1)^2$, and is  biquasiprimitive on vertices while $G^+$ is not quasiprimitive on either bipartite half. Indeed we have $N$ (of order  $2^f(2^{2f}-1)$) normal in $G^+$ and intransitive on each bipartite half (Proposition~\ref{prop-basic}). 
These graphs are quite large, indeed their number of vertices is $2^{2f}(2^{2f}-1)^2/3$ (Proposition~\ref{prop:generalities}).
Infinitely many of them are connected (Proposition~\ref{numberisom}). 
The number of pairwise non-isomorphic connected graphs produced by Construction~\ref{cons} grows exponentially with $f$ (Proposition~\ref{numberisom}); and each connected graph has relatively large girth (at least 10, Proposition~\ref{prop-arc}) and diameter (at least $6f-3$, Proposition~\ref{prop:generalities}).

Note that $G$ is not the full automorphism group of $\Gamma(f,\alpha)$. Moreover overgroups of biquasiprimitive and quasiprimitive groups are not necessarily biquasiprimitive or quasiprimitive respectively. Indeed we have the following:

\begin{theorem}\label{thm:curiosity}
For each connected graph $\Ga=\Gamma(f,\alpha)$ defined in Construction~{\rm\ref{cons}}, with automorphism group $A=\Aut(\Gamma)$ given in Proposition~{\rm\ref{prop:aut}}, $G$ is an index two subgroup of $A$, $\Gamma$ is $(A,3)$-arc-transitive, $A$ is not biquasiprimitive on vertices and $A^+$ is quasiprimitive on each bipartite half.
\end{theorem}

We do not know if there are examples where $G$ is the full automorphism group. 
\begin{question}
Is there a $(G,2)$-arc transitive graph $\Ga$ such that  $G=\Aut(\Ga)$ is biquasiprimitive on vertices but $G^+$ is not quasiprimitive on each bipartite half?  
\end{question}

\section*{Acknowledgements}
The authors thank an anonymous referee for advice that included drawing our attention to the construction in Section \ref{sec:1arc}.

\section{Preliminary graph definitions}
\label{sec:prelim}

We consider simple, undirected graphs $\Ga$, with vertex-set $V\Ga$ and edge-set $E\Ga$. A graph is called {\it cubic} if it is regular of valency $3$.
For a positive integer $s$, an {\it $s$-arc} of a graph is an $(s+1)$-tuple $(v_0,v_1,\ldots, v_s)$ of vertices such that $v_i$ is adjacent to $v_{i-1}$ for all $1\leqslant i\leqslant s$ and $v_{j-1}\neq v_{j+1}$ for all $1\leqslant j\leqslant s-1$. 
The {\it distance} between two vertices $v_1$ and $v_2$, denoted by $\dist_\Gamma(v_1,v_2)$, is the minimal number $s$ such that there exists an $s$-arc between $v_1$ and $v_2$. 
For a connected graph $\Ga$, we define the \textit{diameter of $\Ga$}, denoted $\diam(\Ga)$, as the maximum distance between two vertices of $\Ga$.

We denote a complete graph on $n$ vertices by $K_n$ and a complete bipartite graph with bipartite halves of sizes $n$ and $m$ by $K_{n,m}$. The disjoint union of $m$ copies of $\Sigma$ is denoted by $m\Sigma$. 

Let $\Ga$ be a graph, $G\leqslant \Aut(\Ga)$, and $N\lhd G$. The {\it (normal) quotient graph} $\Ga_N$ is the graph with vertex-set the set of $N$-orbits, such that two $N$-orbits $B_1$ and $B_2$ are adjacent in  $\Ga_N$ if and only if there exist $v\in B_1$ and $w\in B_2$ with  $\{v,w\}\in E\Ga$.

Tables \ref{def:trans} and \ref{def:localtrans} describe some properties $\calP$ that hold for the $G$-action on a connected graph $\Gamma$, where $G\leqslant\Aut(\Gamma)$ and we require that $G$ be transitive on each set in some collection $\calP(\Gamma)$ of sets. For the local variant we require that for each vertex $v$ of $\Gamma$, the stabiliser $G_v$ be transitive on each set in a related collection $\calP(\Gamma,v)$ of sets.  These concepts are sometimes used without reference to a particular group $G$, especially when $G=\Aut(\Ga)$.

\begin{table}
\begin{center}
\begin{tabular}{|l|l|}
\hline
Property   & $\calP(\Ga)=\{\Delta_i|1\leq i\leq s\}$, $\Delta_s\neq\emptyset$\\
\hline 
$(G,s)$-arc transitivity& $\Delta_i$ is the set of $i$-arcs of $\Ga$\\ 
$G$-arc transitivity& $s=1$ and $\Delta_1$ is as in previous line\\
$(G,s)$-distance transitivity& $\Delta_i$ is $\{(v,w)\in V\Gamma\times V\Gamma|\dist_\Gamma(v,w)=i\}$\\ 
$G$-distance transitivity& $s=\diam(\Ga)$ and $\Delta_i$  is as in previous line\\ 
\hline
\end{tabular}
\end{center}
\caption{Properties for $G$-action on a connected graph $\Gamma$}\label{def:trans}
\end{table}

\begin{table}
\begin{center}
\begin{tabular}{|l|l|}
\hline
Local property   & $\calP(\Gamma,v)=\{\Delta_i(v)|1\leq i\leq s\}$, $\Delta_s(v)\neq\emptyset$ for some $v$\\
\hline 
local $(G,s)$-arc transitivity& $\Delta_i(v)$ is the set of $i$-arcs of $\Ga$ with initial vertex $v$\\
local $G$-arc transitivity& $s=1$ and $\Delta_1(v)$ is as previous line\\
local $(G,s)$-distance transitivity& $\Delta_i(v)$ is $\Gamma_i(v):=\{w\in V\Gamma|\dist_\Gamma(v,w)=i\}$\\
local $G$-distance transitivity& $s=\diam(\Ga)$ and $\Delta_i(v)$  is as in previous line\\
\hline
\end{tabular}
\end{center}
\caption{Local properties for $G$-action on a connected graph $\Gamma$}\label{def:localtrans}
\end{table}

Next we describe coset graphs, which will be used to describe our family of graphs, and some of their properties.
\begin{definition}\label{cosetgraph}{\rm
Given a group $G$, a subgroup $H$ and an element $g\in G$ such that $HgH=Hg^{-1}H$, the {\it coset graph} $\Cos(G,H,HgH)$ is the graph with vertices the right cosets of $H$ in $G$, with $Hg_1$ and $Hg_2$ forming an edge if and only if $g_2g_1^{-1}\in HgH$.}
\end{definition}
Note that a coset graph is indeed undirected since $g_2g_1^{-1}\in HgH$ if and only if $g_1g_2^{-1}\in Hg^{-1}H$.

\begin{lemma}\label{lem:cosetgraph} Let $\Ga=\Cos(G,H,HgH)$. Then the following facts hold.
\begin{itemize}
\item[(a)] $\Ga$ has  $|G:H|$ vertices and is regular with valency $|H:H^g\cap H|$.
\item[(b)]  The group $G$ acts by right multiplication on the coset graph with kernel $\cap_{x\in G}H^x$, and $G$ is arc-transitive.
\item[(c)] $\Ga$ is connected if and only if $\la H,g\ra=G$.
\item[(d)] If $\la H,g\ra\leqslant K < G$, then $\Ga=m\Sigma$ where $m=|G:K|$ and $\Sigma=\Cos(K,H,HgH)$. 
\item[(e)] $\Ga$ has $|G:\la H,g\ra|$ connected components, each isomorphic to  $\Cos(\la H,g\ra,H,HgH)$. 
\item[(f)] For $\eta\in\N_{\Aut\, G}(H)$, the map $\bar\eta:Hx\mapsto Hx^\eta$ is a permutation of $V\Ga$ and induces an isomorphism from $\Ga$ to $\Cos(G,H,Hg^\eta H)$.
\end{itemize}
 \end{lemma}
\proof Statements (a) to (c) can be found in \cite{Lor}.

Assume $\la H,g\ra\leqslant K < G$.
By Theorem 4(i,iii) of \cite{Lor}, there is no edge of $\Ga$ between vertices (that is, $H$-cosets) lying in distinct $K$-cosets.
On the other hand, by the last paragraph of the proof of that same theorem, for all $K$-cosets $Kx$, the graph induced on the $H$-cosets contained in $Kx$ is isomorphic to $\Sigma=\Cos(K,H,HgH)$.
Hence (d) holds. Statement (e) follows from (d) (taking $K=\la H,g\ra$) and (c).

Let  $\eta\in\N_{\Aut\, G}(H)$ and $\Sigma=\Cos(G,H,Hg^\eta H)$. Then $\eta$ maps $H$-cosets to $H$-cosets and so induces the permutation $\bar\eta:V\Ga \rightarrow V\Ga:Hx\mapsto Hx^\eta$ of $V\Ga=V\Sigma$. 
Let $\{Hx,Hy\}$ be an edge of $\Ga$, that is, $yx^{-1}\in HgH$. Now $y^\eta(x^\eta)^{-1}=(yx^{-1})^\eta\in (HgH)^\eta$. Since $\eta$ normalises $H$, we have $(HgH)^\eta=Hg^\eta H$, and so $\{Hx^\eta,Hy^\eta\}$ is an edge of $\Sigma$. Conversely, let $\{Hx^\eta,Hy^\eta\}$ be an edge of  $\Sigma$, so that $y^\eta(x^\eta)^{-1}=(yx^{-1})^\eta\in Hg^{\eta}H$. Then $yx^{-1}\in (Hg^\eta H)^{\eta^{-1}}$, and since $\eta$ normalises $H$, $(Hg^\eta H)^{\eta^{-1}}=HgH$. Therefore $\bar\eta$ sends the edge-set of $\Ga$ to the edge-set of  $\Sigma$ and (f) holds.
\qed

\section{$1$-arc-transitive examples}
\label{sec:1arc}

In this section we construct an infinite family of $G$-arc-transitive graphs such that $G$ is biquasiprimitive on vertices but $G^+$ is not quasiprimitive on each bipartite half.

\begin{construction}\label{s=1}
\textrm{
 Let $G=\mathbb{Z}_p\times \mathbb{Z}_p\times \mathbb{Z}_2$ where $p \equiv 1 \pmod 3$ is a prime. Let $a$ be an element of multiplicative order $3$ in $\mathbb{Z}_p$.
We define a graph $\Sigma$ with vertex-set $G$ and edges of the form
$$\begin{array}{ll}
 \{(x,y,0), &(x+1,y+1,1)\}, \\
  \{(x,y,0), &(x+a,y+a^2,1)\}, \\
 \{(x,y,0), &(x+a^2,y+a,1)\}
\end{array}$$
This yields an undirected bipartite graph with bipartite halves $\Delta_1=\{(x,y,0)|x,y\in\mathbb{Z}_p\}$ and $\Delta_2=\{(x,y,1)|x,y\in\mathbb{Z}_p\}$.
}%

\textrm{
Some automorphisms of $\Sigma$ are:
\begin{itemize}
 \item $t_{u,v}:(x,y,\epsilon)\mapsto (x+u,y+v,\epsilon)\in \Aut(\Sigma)$, we denote $\{t_{u,v}|u,v\in \mathbb{Z}_p\}$ by $N\cong \mathbb{Z}_p^2 $;
\item $\tau:(x,y,\epsilon)\mapsto (ax,a^2y,\epsilon) \in \Aut(\Sigma)$;
\item $\sigma: (x,y,\epsilon)\mapsto (y,x,\epsilon) \in \Aut(\Sigma)$;
\item $\nu: (x,y,\epsilon)\mapsto (-x,-y,1-\epsilon) \in \Aut(\Sigma)$.  
\end{itemize}
}%
\end{construction}

\begin{proposition}
\label{prn:s=1}
Let  $G=N\rtimes\langle \tau,\sigma\nu\rangle\cong \mathbb{Z}_p^2\rtimes S_3$.
Then $G$ is biquasiprimitive on $V\Sigma$ but $G^+$ is not quasiprimitive on each bipartite half and $\Sigma$ is $(G,1)$-arc transitive but not $(G,2)$-arc transitive. The full automorphism group of $\Sigma$ is  $A=N\rtimes\langle \tau,\sigma,\nu\rangle\cong \mathbb{Z}_p^2\rtimes (S_3\times\mathbb{Z}_2)$. Then $\Sigma$ is $(A,2)$-arc transitive, $A$ is biquasiprimitive on $V\Sigma$ and $A^+$ is quasiprimitive on the bipartite halves.
\end{proposition}
\begin{proof}
The group $N$ clearly acts transitively on each bipartite half and $\sigma\nu$ switches $\Delta_1$ and $\Delta_2$, so $G$ is transitive on $V\Sigma$. Moreover, since no nontrivial element of  $\langle \tau,\sigma\nu\rangle$ centralises $N$, it follows that $N$ is the unique minimal normal subgroup of $G$ and so $G$ is biquasiprimitive on vertices. Now $G^+=N\rtimes\langle \tau\rangle$ has $\{t_{u,0}|u\in \mathbb{Z}_p \}\cong\mathbb{Z}_p$ as a normal subgroup that is intransitive on $\Delta_1$ and $\Delta_2$. Thus $G^+$ is not quasiprimitive on each bipartite half. Finally, $G_{(0,0,0)}=\langle\tau\rangle$, which acts regularly on the set of three neighbours of $(0,0,0)$ and so $\Sigma$ is $(G,1)$-arc transitive but not $(G,2)$-arc transitive.

Let $A=N\rtimes\langle \tau,\sigma,\nu\rangle$. Then $N$ is also a unique minimal normal subgroup of $A$ and of $A^+=N\rtimes\langle \tau,\sigma\rangle$. Thus $A$ is biquasiprimitive on vertices and $A^+$ is quasiprimitive on each bipartite half. Moreover,  $A_{(0,0,0)}=\langle \tau,\sigma\rangle\cong S_3$ acts 2-transitively on the set of three neighbours of $(0,0,0)$ and so $\Sigma$ is $(A,2)$-arc transitive. 

Let $X$ be the full automorphism group of $\Sigma$. Since $A$ is vertex-transitive we have $X=AX_{\alpha}$ (where $\alpha\in V\Sigma$) and so $|X_{\alpha}|$ divides 48 \cite{Tutte1,Tutte2}. Since $|A_{\alpha}|=6$ it follows that $|X:A|$ divides 8. Consider the action of $X$ on the set of right cosets of $A$. If $A$ is core-free in $X$ it follows that $X\leqslant S_8$, contradicting $p^2$ dividing $|A|$ and $p\geq 7$. Thus $A$ contains a normal subgroup $M$ of $X$. Since $N$ is the unique minimal normal subgroup of $A$ it follows that $N\leqslant M$. However, $N$ is the unique Sylow $p$-subgroup of $A$ and hence of $M$, and so $N\norml X$. Hence $X$ has a normal subgroup that acts regularly on each bipartite half and so by \cite[Lemma 2.4]{li}, $X_{\alpha}$ acts faithfully on $\Sigma(\alpha)$. Thus $X_{\alpha}=A_{\alpha}=S_3$ and hence $X=A$.
\end{proof}

\section{Finite fields}
This section contains facts about finite fields that we need later. We denote a field of order $q$ by $\GF(q)$.

\begin{definition}{\rm Let $x$ be an element of a field $F$. 
The {\it subfield generated by} $x$ is the unique smallest subfield containing $x$.
The element $x$ is called a {\it generator} of $F$ if the subfield generated by $x$ is $F$, in other words, if $x$ is not contained in any proper subfield of $F$. 
}\end{definition}

\begin{lemma}\label{lem:gen}
Let $f$ be an integer and let $\alpha\in \GF(2^f)$.  
The subfield generated by $\alpha$ is $\GF(2^e)$  if and only if the  order of $\alpha$ divides $2^e-1$ but does not divide $2^s-1$ for any proper divisor $s$ of $e$.
In particular, $\alpha$ is a generator of $\GF(2^f)$ if and only if the  order of $\alpha$ does not divide $2^e-1$ for any proper divisor $e$ of $f$.
\end{lemma}
\proof
Since the multiplicative group of $\GF(2^f)$ is cyclic of order $2^f-1$, it follows that the multiplicative group of the subfield $\GF(2^e)$ of $\GF(2^f)$ is precisely the subgroup of order $2^e-1$, with $e$ dividing $f$. This subgroup is unique, since there is a unique subgroup of each order in a cyclic group. Thus the order of $\alpha$ divides  $2^e-1$ if and only if $\alpha\in \GF(2^e)$. The result follows.
\qed

\begin{lemma}\label{lem:alpha+1} Let $f$ be an integer, $f\geq 2$, and let $\alpha$ be a generator of $\GF(2^f)$. Then
\begin{itemize}
 \item[(1)]   $\alpha^{2^i}\neq \alpha+1$ for all positive integers $i<f$ except possibly $i=f/2$ (with $f$ even), and
\item[(2)] $\alpha^{2^i}\neq \alpha$ for all positive integers $i<f$.
\end{itemize}
 \end{lemma}
\proof
  Suppose $\alpha^{2^i}= \alpha+1$ for some integer $i<f$.
Then since $\GF(2^f)$ has characteristic $2$, we have  $\alpha^{2^{2i}}=(\alpha^{2^{i}})^{2^i}=(\alpha+1)^{2^i}=\alpha^{2^{i}}+1=\alpha$, so  $\alpha^{2^{2i}-1}=1$. 
Since $0\neq \alpha\in \GF(2^f)$, we also have that $\alpha^{2^{f}-1}=1$. Hence the order of $\alpha$ divides $\gcd(2^{2i}-1,2^f-1)=2^{\gcd(2i,f)}-1$. Since $\gcd(2i,f)$ divides $f$ and $\alpha$ is a generator, Lemma \ref{lem:gen} implies that $\gcd(2i,f)=f$, that is, $f$ divides $2i$. Since $f>i$, this implies that $f$ is even and $i=f/2$. This proves (1).

Suppose  $\alpha^{2^i}= \alpha$ for some positive integer $i<f$. Then $\alpha^{2^{i}-1}=1$. Hence the order of $\alpha$ divides $\gcd(2^{i}-1,2^f-1)=2^{\gcd(i,f)}-1$. Since $\gcd(i,f)$ is a divisor of $f$ and $\alpha$ is a generator, Lemma \ref{lem:gen} implies that $\gcd(i,f)=f$, that is, $f$ divides $i$, contradicting $f>i$. This proves (2).
\qed 

\begin{lemma}\label{ff}
Let $f$ be an integer, $f\geq 3$. Then the number of generators of $\GF(2^{f})$ is strictly greater than $2^{f-1}$.
\end{lemma}
\proof
For $f=3$, all elements of  $\GF(2^{3})\setminus\{0,1\}$ are generators, hence there are $6$ generators and the claim holds.
Assume $f\geq 4$.
Let $f=\Pi_{i=1}^k p_i^{e_i}$, where the $p_i$ are distinct primes and each $e_i\geq 1$. Let $f_i=f/p_i$.
Then all elements which are not generators are in one of the subfields $\GF(2^{f_i})$. Hence the number of generators is $2^f-|\cup_{i=1}^k \GF(2^{f_i})|$.
We have $|\cup_{i=1}^k \GF(2^{f_i})|\leq 1+\Sigma_{i=1}^k (2^{f_i}-1)$ since $0$ is in all fields. Since $f_i\leq f/2$ for all $i$, we have $|\cup_{i=1}^k \GF(2^{f_i})|\leq 1+ k(2^{f/2}-1)\leq  k2^{f/2}$.
Since $f\geq \Pi_{i=1}^k p_i\geq 2^k$, we have $k\leq \log_2(f)$, and so  $|\cup_{i=1}^k \GF(2^{f_i})|\leq  \log_2(f)2^{f/2}$.
It is easy to check that, for $f\geq 4$,  $\log_2(f)\leq 2^{f/2-1}$, and so $\log_2(f)2^{f/2}\leq 2^{f-1}$.
We can now conclude that the number of generators is at least $2^f-2^{f-1}=2^{f-1}$. 

Suppose we get equality. Then we have equality in all our inequalities. In particular $1+ k(2^{f/2}-1)=  k2^{f/2}$, and so $k=1$, and  $k= \log_2(f)$, so $f=2^k$. Thus $f=2$, a contradiction. 
Therefore the number of generators is greater than $2^{f-1}$.
\qed

\begin{lemma}\label{ffeven}
Let $\ell$ be an integer, $\ell\geq2$. Then the number of generators of $\GF(2^{2\ell})$ which do not satisfy the equation $x^{2^{\ell}}=x+1$ is strictly greater than $2^\ell(2^{\ell-1}-1)$ 
\end{lemma}
\proof
By Lemma \ref{ff}, $\GF(2^{2\ell})$ contains more than  $2^{2\ell-1}$ generators.
Since the equation $x^{2^{\ell}}=x+1$ has degree $2^{\ell}$, it has at most $2^{\ell}$ solutions. Hence the number of generators not satisfying the equation is greater than $2^{2\ell-1}- 2^\ell = 2^\ell(2^{\ell-1}-1)$.
\qed

\section{The group $\PSL(2,2^f)$}\label{sec:group}

The elements of a group $\PSL(2,q)$ may be viewed as permutations of $X:=\GF(q)\cup\{\infty\}$. More precisely $t_{a,b,c,d}$ is the element:
\begin{equation}\label{gabc}
{t_{a,b,c,d}}:\ x\ \mapsto  \frac{ax+b}{cx+d}\quad \mbox{for all $x \in X$}
\end{equation}
where $a,b,c,d\in \GF(q)$ are such that $ad-bc$ is a nonzero square of $\GF(q)$. We adopt the convention that $\infty$ is mapped by ${t_{a,b,c,d}}$ onto $ac^{-1}$ and that an element of $\GF(q)$ divided by $0$ is $\infty$.
For $q=2^f$, all nonzero elements of $\GF(q)$ are squares, and the automorphism group of $\PSL(2,q)$ is $\PGammaL(2,q) = \la \PSL(2,q),\tau\ra$, where
\begin{equation}\label{tau}
\tau:t_{a,b,c,d} \mapsto t_{a^2,b^2,c^2,d^2} \quad\mbox{for each 
$t_{a,b,c,d}\in \PSL(2,q)$.}
\end{equation}

In this paper we will take $T=\PSL(2,2^f)$ for some $f\geq 1$. 
For each subfield $\GF(2^e)$ of $\GF(2^f)$, we identify $\PSL(2,2^e)$ with the subgroup of $T$ of those $t_{a,b,c,d}$ with all of $a,b,c,d\in \GF(2^e)$.
In our construction, we will use the following notation for elements of $H=\PSL(2,2)\leqslant T$. 
\begin{equation}\label{notations}
a=t_{1,1,1,0}:\: x\mapsto 1+\frac{1}{x}, \quad b=t_{1,1,0,1}:\: x\mapsto x+1.
\end{equation}
Note that $a$ has order 3, $b$ has order 2, and $H=\langle a\rangle\rtimes\la b\ra\cong S_3$.
For $\alpha\in \GF(2^f)$, we will also need the following elements of $T$:
\begin{equation}\label{notations2}
 u_\alpha=t_{1,\alpha,0,1}: \: x\mapsto x+\alpha, \quad c_\alpha=a^{u_\alpha}=t_{\alpha+1,\alpha^2+\alpha+1,1,\alpha}.
\end{equation}
Let $P$ be the Sylow 2-subgroup of $T$ containing the involution $b=u_1$, that is,  $P=\{u_\alpha|\alpha\in \GF(2^f)\}$.
Then $\N_T(P)\cong\AGL(1,2^f)$ is the set of permutations $t_{r,s,0,1}:\: x\mapsto rx+s$ with $r\neq 0$.

\begin{lemma}\label{lem:groupprop}
Let  $\alpha\in \GF(2^f)$.
Using the notation introduced above, the following facts hold.
 \begin{itemize}
\item[(a)] $\C_T(b)=P$. In particular,  $u_\alpha b=bu_\alpha=u_{\alpha+1}$  and $\C_H(b)=\langle b \rangle$.
\item[(b)]  For $\alpha\neq 0$, the element $z_\alpha:=t_{\alpha^{-1},0,0,1}\in\N_T(P)$. Moreover 
$u_\alpha=b^{z_\alpha}$ and the order of $z_\alpha$ is equal to the multiplicative order of $\alpha$.
\item[(c)] $c_\alpha^{\tau^i}=c_\alpha^{-1}$ if and only if $\alpha^{2^i}=\alpha+1$.
\item[(d)] $\N_T(H)=H$.
\item[(e)] If the subfield generated by $\alpha$ is $\GF(2^e)$, then $\la H,\ua\ra=\PSL(2,2^e)$. 
\end{itemize}
\end{lemma}
\proof
(a) The centraliser of $b$ in $T$ is easily computed. Since $\ua\in P$, it then commutes with $b$, and $b\ua=u_{\alpha+1}$. Also $\C_H(b)=\C_T(b)\cap H=\langle b \rangle$. \\
(b) A calculation shows that $u_y^{z_\alpha}=u_{\alpha y}\in P$, and so $z_\alpha\in \N_T(P)$. Also $u_{\alpha}=u_1^{z_\alpha}=b^{z_\alpha}$.
Since $z_\alpha^j=t_{\alpha^{-j},0,0,1}$ the rest of the statement follows. \\
(c) This is a simple calculation left to the reader.\\
(d) Let $D=\N_T(\langle a\rangle)$. Now $D$ is a dihedral group  $D_{2(2^f\pm 1)}$, see \cite[Section 260]{Dickson}. Since $\langle a\rangle\cong C_3$ is characteristic in $H\cong S_3$,  $\N_T(H)\leqslant \N_T(\langle a\rangle)=D$, and so $\N_T(H)=\N_D(H)$. Since an $S_3$ subgroup in a dihedral group $D_{2n}$, $n$ odd, is self-normalising, we have that $\N_D(H)=H$. Thus  $\N_T(H)=H$.\\ 
(e) Suppose the subfield generated by $\alpha$ is $\GF(2^e)$. 
If $e=1$, then $\alpha=0$ or $1$, $u_\alpha\in H$ and  $\la H,\ua\ra=H=\PSL(2,2)$.
Assume now $e\geq 2$.
Since all the subscripts of $\ua=t_{1,\alpha,0,1}$ are in $\GF(2^e)$, we obviously have $\la H,\ua\ra\leqslant\PSL(2,2^e)$.
Suppose that $\la H,\ua\ra\leqslant M$, where $M$ is a maximal subgroup of $\PSL(2,2^e)$. Since $\la H,\ua\ra$ contains a subgroup isomorphic to $S_3$, $M$ cannot be isomorphic to $\AGL(1,2^e)$ (for $e$ even, 3 divides $|\AGL(1,2^e)|$ but no involution in $\AGL(1,2^e)$ inverts an element of order 3). Also since $\la H,\ua\ra$ contains subgroups which are isomorphic to $C_2^2$, $M$ cannot be isomorphic to $D_{2(2^e\pm 1)}$. It follows from the list of maximal subgroups of $\PSL(2,2^e)$ (see \cite[Section 260]{Dickson}) that  $M\cong \PSL(2,2^s)$ for some proper divisor  $s$  of $e$.
Since $b,\ua\in M$ and commute, they lie in the same Sylow 2-subgroup $S$ of $M$, so there exists $d\in M$ such that $b^d=\ua$. Hence $b^d=\ua=b^{z_\alpha}$ (by Part (b)), and so $dz_\alpha^{-1}$ centralises $b$. Since $\C_T(b)=P$ by (a), we obtain that $d\in Pz_\alpha$. Since $z_\alpha\in\N_T(P)$ has order $n:=|\alpha|$, it follows that $d$ has order divisible by $n$. 
 Moreover, $d$ must be in $\N_M(S)\cong \AGL(1,2^s)$, and so the order of $d$ divides $2^{s}-1$. Thus $n$ divides  $2^{s}-1$, a contradiction to Lemma \ref{lem:gen}.
Thus, $\langle H,\ua\rangle=\PSL(2,2^e)$.
\qed

\section{The family of graphs}
Let $f$ be a positive integer, and let $T$, $H$, $a$, $b$, $\alpha$, $z_\alpha$ (for $\alpha\neq 0$), $\ua$, and $c_\alpha$ be as in Section \ref{sec:group}.

\begin{construction}\label{cons}
{\rm 
Let $G=T^2\rtimes\la\pi\ra$, where $\pi\in\Aut(T^2)$ is such that $(x,y)^\pi=(y,x)$, for all elements $(x,y)\in T^2$.
Let $L=\la(a,a),(b,b)\ra<T^2$, and
\begin{equation}\label{defg}
 g_\alpha =(u_\alpha,b u_{\alpha})\pi=(u_\alpha,u_{\alpha} b)\pi=(t_{1,\alpha,0,1}, t_{1,\alpha+1,0,1})\pi.
\end{equation}
By Lemma \ref{lem:(a,a)}(c) below, $\ga^{-1}=\ga (b,b)$. Thus $L\ga^{-1}L=L\ga (b,b)L=L\ga L$.
Define $\Ga=\Ga(f,\alpha)=\Cos(G,L,L\ga L)$.
}
\end{construction}

We shall need information about the following subgroups:
\begin{equation}\label{notations3}
X_\alpha=\la L,\ga\ra,\quad N_\alpha=\la L,(c_\alpha^{-1},c_\alpha)\ra.
\end{equation}

\begin{lemma}\label{lem:(a,a)}
The following facts hold.
\begin{itemize}
\item[(a)] $|G|=2^{2f+1}(2^{2f}-1)^2$.
\item[(b)]  $(a,a)^{g_\alpha}=(c_\alpha^{-1},c_\alpha)$, where $c_\alpha$ is as in  {\rm (\ref{notations2})} and has order $3$. Thus $N_\alpha\leqslant X_\alpha$.
\item[(c)] $\ga^{-1}=\ga (b,b)$ and $(b,b)^{g_\alpha}=(b,b)$. 
\item[(d)] For $f\ge 2$ and $\alpha$ a generator of $\GF(2^f)$, either $N_\alpha= T^2$ or $N_\alpha= \{(t,t^\nu)|t\in T\}\cong T$ for some $\nu\in \Aut(T)$.
\end{itemize}
\end{lemma}
\proof
(a) follows from the fact that $|G|=2|T|^2$.\\
(b) We have $(a,a)^{g_\alpha}=(a^{u_\alpha},(a^b)^{u_\alpha})^\pi=
(c_\alpha,c_\alpha^{-1})^\pi=(c_\alpha^{-1},c_\alpha)$, by (\ref{notations2}), and hence $N_\alpha\leqslant X_\alpha$. Since $c_\alpha$ is conjugate to $a$, it has order 3. \\
(c) We have $\ga^2 (b,b)=(u_\alpha,b u_{\alpha})\pi(u_\alpha,b u_{\alpha})\pi(b,b)=(u_\alpha,b u_{\alpha})(b u_\alpha, u_{\alpha})(b,b)=(1,1)$ since  $u_\alpha b=bu_\alpha$ by Lemma \ref{lem:groupprop}(a). Thus $\ga^{-1}=\ga (b,b)$.
We also have $(b,b)^{g_\alpha}=(b^\ua,b^{\ua b})^\pi=(b,b)^\pi=(b,b)$, using Lemma \ref{lem:groupprop}(a) for the second equality.\\
(d) The projections of $N_\alpha$ onto each of 
the two coordinates are equal to $\la a,b,c_\alpha\ra$.
Since $\ua b=b\ua$, the subgroup $\la a,b,c_\alpha\ra$ of $T$ is normalised by each of $a,b$ and $\ua$.
Hence $\la a,b,c_\alpha\ra\lhd\la a,b,\ua\ra$, and  $\la a,b,\ua\ra=T$ by Lemma \ref{lem:groupprop}(e). Thus $\la a,b,c_\alpha\ra=T$ since $T$ is simple, and so $N_\alpha= T^2$ or $N_\alpha\cong T$. In the latter case, $N_\alpha$ is a diagonal subgroup of $T^2$ and hence $N_\alpha= \{(t,t^\nu)|t\in T\}\cong T$ for some $\nu\in \Aut(T)$.
\qed

We first describe some general properties of the graphs $\Ga(f,\alpha)$.

\begin{proposition}\label{prop:generalities}
Let $f\geqslant 1$ be an integer and $\alpha$ be an element of $\GF(2^f)$.
Let $\Ga=\Ga(f,\alpha)$, $G$, $T$, $L$, $\pi$ be as in Construction \ref{cons}.
Then  $\Ga$ is bipartite, cubic, and, if $\Ga$ is connected, then it has diameter at least $6f-3$. Moreover, $G^+=T^2$, $G\leqslant\Aut(\Ga)$ and $|V\Ga|=2^{2f}(2^{2f}-1)^2/3$. 
\end{proposition}
\proof
By Lemma \ref{lem:(a,a)}(b),  $(a,a)^{g_\alpha}=(c_\alpha^{-1},c_\alpha)$, which is not in $L$ since $c_\alpha\neq c_\alpha^{-1}$,  and, by Lemma \ref{lem:(a,a)}(c), $(b,b)^{g_\alpha}=(b,b)$. Thus the intersection $L^{g_\alpha}\cap L=\la (b,b)\ra\cong C_2$, and so the graph $\Ga$ has valency $|L:L^\ga\cap L|=3$ (hence is cubic) by Lemma \ref{lem:cosetgraph}(a). 
Moreover, $T^2$ has two orbits on the cosets of $L$, and since $T^2\cap L\ga L=\emptyset$, no vertices in the same orbit are adjacent. Hence $\Ga$ is bipartite. Since $T^2$ is an index 2 subgroup of $G$ and its orbits are the two bipartite halves, the even subgroup $G^+$ is precisely $T^2$.  
The number of vertices of $\Ga$ is $|G|/|L|=2^{2f}(2^{2f}-1)^2/3$, with each bipartite half of size $2^{2f-1}(2^{2f}-1)^2/3$. 

Suppose $\Ga$ is connected and let $d=\diam(\Ga)$. We have $|\Ga_1(L)|=3$ and $|\Ga_i(L)|$ is at most $2|\Ga_{i-1}(L)|$ for $2\leqslant i\leqslant d$. Hence the number of vertices of $\Ga$ is at most $1+3+3.2+\ldots+3.2^{d-1}=1+3(2^d-1)$.
Therefore $2^{2f}(2^{2f}-1)^2/3 \leqslant 1+3(2^d-1)$, or equivalently  $2^{2f}(2^{2f}-1)^2/9+2/3 \leqslant 2^d$, which implies $2^{2f}(2^{2f}-1)^2/9 < 2^d$. Thus $(2^{2f}-1)/3 < 2^{\frac{d}{2}-f}$. 
Now for all $f\geq 1$, we have  $(2^{2f}-1)/3\geq 2^{2f}/4=2^{2f-2}$, and so  $2^{2f-2}<2^{\frac{d}{2}-f}$. Therefore  $2f-2<\frac{d}{2}-f$ and $d>6f-4$.
Since $\cap_{x\in G}L^x$ is trivial, it follows from Lemma \ref{lem:cosetgraph}(b) that $G$ acts faithfully on $\Ga$, and hence $G\leqslant \Aut(\Ga)$.
\qed

Note that the bound on the diameter is not tight. For example, for $f=3$ a MAGMA \cite{magma} computation shows that $\Ga(3,\alpha)$ has diameter 21 for $\alpha$ a generator of $\GF(8)$ (we will see in Corollary \ref{cor:f=3} that the graph is connected in this case).

\section{Equality and connectivity}

We first have a lemma determining when graphs obtained by Construction \ref{cons} have the same edge-set. 

\begin{proposition}\label{prop:samegraph}
Let $f\ge 1$. Let $\alpha,\beta$ be  elements of $\GF(2^f)$. 
Then $\Ga(f,\alpha)=\Ga(f,\beta)$ if and only if  $\beta\in\{\alpha,\alpha+1\}$. 
\end{proposition}
\proof
Suppose that $\Ga(f,\alpha)=\Ga(f,\beta)$. Then the double cosets $Lg_\alpha L$ and $Lg_\beta L$ coincide, and so $g_\beta\in Lg_\alpha L$.
Since $\pi$ centralises $L$, this implies, using (\ref{defg}), that  
$(u_\beta,b u_\beta)=(h_1,h_1)(u_\alpha,b u_\alpha)(h_2,h_2)$ for some $h_1,h_2\in H$.
Thus $h_1b \ua h_2=b u_\beta=bh_1\ua h_2$, and so $h_1$  commutes with $b$.
Since $b$ centralises $P$ by Lemma \ref{lem:groupprop}(a) and $\ua,u_\beta\in P$, we also have $h_1 \ua b h_2= u_\beta b=h_1\ua h_2 b$, and so $h_2$ also  commutes with $b$.
Hence $h_1,h_2\in \C_H(b)=\langle b \rangle$ by Lemma \ref{lem:groupprop}(a). If $h_1=h_2$, then $\alpha=\beta$, and if $h_2=h_1 b$ then $\beta=\alpha +1$. 

Conversely, if $\beta=\alpha+1$, then $g_\beta=(u_\beta,u_\beta b)\pi=(u_\alpha b,u_\alpha)\pi=\ga (b,b)$, and so  $Lg_\alpha L=Lg_\beta L$. Thus $\Ga(f,\alpha)=\Ga(f,\beta)$.
\qed

For $f=1$ Construction \ref{cons} yields only one graph.
\begin{lemma}\label{lem:f=1}
 $\Ga(1,0)=\Ga(1,1)=2 \,K_{3,3}$
\end{lemma}
\proof
Here $T=H$, and by Proposition  \ref{prop:samegraph}, $\Ga(1,0)=\Ga(1,1)$ so we may assume $\alpha=0$. Thus $u_\alpha=1$ and $\ga=(1,b)\pi$.
It can be computed that $\la L,\ga \ra=\{(x,y)|x^{-1}y\in \la a\ra\}\cup\{(x,yb)\pi|x^{-1}y\in \la a\ra\}$ has index 2 in $G$. Therefore by Lemma \ref{lem:cosetgraph}(e), 
 $\Ga(1,0)$ has $2$ connected components. Each must be bipartite and have valency 3 by Proposition \ref{prop:generalities}, hence the conclusion. 
\qed

The next two general results allow us to determine the connected components of $\Ga(f,\alpha)$.

\begin{lemma}\label{lem:gensubfield}
Let $\alpha$ be an element of $\GF(2^f)$ and let $\GF(2^e)$ be the subfield generated by $\alpha$.  
Then $\Ga(f,\alpha)\cong m\Ga(e,\alpha)$, where $m=|T:\PSL(2,2^e)|^2$.
\end{lemma}
\proof 
Let $K=\PSL(2,2^e)^2\rtimes\la\pi\ra$ viewed as a subgroup of $G$. Then $\ga\in K$ and  $L\leqslant K$, and so $\la L,\ga\ra\leqslant K$.
By Lemma \ref{lem:cosetgraph}(d), $\Ga(f,\alpha)=m\Sigma$ where $m=|G:K|$ and $\Sigma=\Cos(K,L,L\ga L)$.  Finally, $m=|G:K|=2|T|^2/(2|\PSL(2,2^e)|^2)=|T:\PSL(2,2^e)|^2$.
\qed

\begin{proposition}\label{prop:connected}
Let $f\ge2$ and $\alpha\in \GF(2^f)$ be a generator. 
\begin{itemize}
 \item[(a)] If $f$ is odd, or if $f$ is even and $\alpha^{2^{(f/2)}}\neq\alpha+1$, then $\Ga(f,\alpha)$ is connected.
\item[(b)] If $f$ is even and $\alpha^{2^{(f/2)}}=\alpha+1$, then $\Ga(f,\alpha)$ has $|T|$ connected components, each containing $|T|/3$ vertices and isomorphic to 
$\Cos(\la T,\nu\ra,H,H\ua\nu H)$ where $H=\PSL(2,2)$ and $\nu=\tau^{(f/2)}$.
\end{itemize}
\end{proposition}
\proof
We set $X_\alpha=\la L,\ga\ra$ and $N_\alpha=\la L,(c_\alpha^{-1},c_\alpha)\ra$ as in (\ref{notations3}). By Lemma \ref{lem:cosetgraph}(e),
 the number of connected components of $\Ga(f,\alpha)$ is $|G:X_\alpha|$ and all connected components are isomorphic to $\Cos(X_\alpha,L,L\ga L)$.

We have $\alpha\notin\{0,1\}$, since $\alpha$ is a generator and $f\neq 1$.

By Lemma \ref{lem:(a,a)}(b),   $N_\alpha\leqslant X_\alpha$, and by Lemma \ref{lem:(a,a)}(d), either $N_\alpha= T^2$ or $N_\alpha=\{(t,t^\nu)|t\in T\}$ for some $\nu\in\Aut(T)$.  
In the latter case, since $N_\alpha$ contains $(a,a)$, $(b,b)$ and $(c_\alpha^{-1},c_\alpha)$,  $\nu$ must be in  $\C_{\Aut(T)}(\langle a,b\rangle)$ and must satisfy  $c_\alpha^\nu=c_\alpha^{-1}$. 
Since $\la a,b\ra=\PSL(2,2)$, we have $\C_{\Aut(T)}(\langle a,b\rangle)=\C_{\Aut(T)}(\PSL(2,2))=\Aut(\GF(2^f))=\la\tau\ra\cong C_f$, where $\tau$ is the Frobenius automorphism described in (\ref{tau}).
.

Assume $f$ is odd, or $f$ is even and $\alpha^{2^{(f/2)}}\neq\alpha+1$. Then by Lemma \ref{lem:alpha+1}(1), $\alpha^{2^i}\neq \alpha+1$ for all $i<f$, and so by Lemma \ref{lem:groupprop}(c), $c_\alpha^{\tau^i}\neq c_\alpha^{-1}$ for all $i<f$.  Hence there is no $\nu\in \C_{\Aut(T)}(\langle a,b\rangle)$ satisfying $c_\alpha^\nu=c_\alpha^{-1}$. Thus $N_\alpha= T^2$, and so $X_\alpha=G$ since $\ga\notin T^2$. Thus $\Ga(f,\alpha)$ is connected and (a) holds.

Now assume $f$ is even and $\alpha^{2^{i}}=\alpha+1$, where $i=f/2$. 
Let $\nu:=\tau^i$. By Lemma \ref{lem:groupprop}(c), $\nu\in \C_{\Aut(T)}(\langle a,b\rangle)$ and satisfies $c_\alpha^\nu=c_\alpha^{-1}$, and so $N_\alpha=\{(t,t^\nu)|t\in T\}\cong T$.  Notice $\nu$ is an involution. 
We have $N_\alpha\leqslant X_\alpha$, and so $\la N_\alpha,\ga\ra\leqslant \la X_\alpha,\ga\ra=X_\alpha$. On the other hand, $X_\alpha=\la (a,a),(b,b),\ga\ra\leqslant \la (a,a),(b,b),(c_\alpha^{-1},c_\alpha),\ga\ra=\la N_\alpha,\ga\ra$. Thus $X_\alpha=\la N_\alpha,\ga\ra$. 
Notice that $\ua^\nu=t_{1,\alpha^{2^i},0,1}=u_{\alpha+1}=\ua b$, and so $\ga=(\ua,\ua^\nu)\pi$. Therefore $\la N_\alpha,\ga\ra=\la N_\alpha,\pi\ra=N_\alpha\rtimes\la\pi\ra$.
Hence, $|X_\alpha|=2|N_\alpha|=2|T|$. Moreover, $X_\alpha=\{(t,t^\nu)\pi^\epsilon|t\in T, \epsilon\in\{0,1\}\}$. Also the number of connected components is $|G:X_\alpha|=|T|$ by Lemma \ref{lem:cosetgraph}(e).

We now prove that $X_\alpha$ is isomorphic to $\la T,\nu\ra$.
We define $$\phi:X_\alpha\rightarrow \la T,\nu\ra: (t,t^\nu)\pi^\epsilon\mapsto t \nu^\epsilon.$$

We first show that $\phi$ is a homomorphism, that is, that $\phi((t_1,t_1^\nu)\pi^{\epsilon_1}(t_2,t_2^\nu)\pi^{\epsilon_2})=\phi((t_1,t_1^\nu)\pi^{\epsilon_1})\phi((t_2,t_2^\nu)\pi^{\epsilon_2})$. This clearly holds for $\epsilon_1=0$. We now prove the case $\epsilon_1=1$. 
\begin{eqnarray*}
\phi((t_1,t_1^\nu)\pi(t_2,t_2^\nu)\pi^{\epsilon_2})&=& \phi((t_1,t_1^\nu)(t_2^\nu,t_2)\pi\pi^{\epsilon_2})\\
&=&\phi((t_1 t_2^\nu,t_1^\nu t_2)\pi^{1-\epsilon_2})\\
&=&t_1 t_2^\nu\nu^{1-\epsilon_2}\\
&=&t_1 \nu t_2\nu\nu^{1-\epsilon_2}\\
&=&(t_1 \nu) (t_2 \nu^{\epsilon_2})\\
&=&\phi((t_1,t_1^\nu)\pi)\phi((t_2,t_2^\nu)\pi^{\epsilon_2}).
\end{eqnarray*}
Thus $\phi$ is a homomorphism. Clearly $\Ker\phi=1$, and $|X_\alpha|=|\la T,\nu\ra|=2|T|$, and so $\phi$ is a bijection. Therefore $\phi$ is an isomorphism.

Notice that $\phi(L)=\la a,b\ra=H$ and $\phi(\ga)=\ua\nu$.

By Lemma \ref{lem:cosetgraph}(e), each connected component of $\Ga(f,\alpha)$ is isomorphic to $\Cos(X_\alpha,L,L\ga L)$, and $\phi$ induces a graph isomorphism  
$\Cos(X_\alpha,L,L\ga L)\cong\Cos(\la T,\nu\ra,H,H\ua\nu H)$. Thus (b) holds.
\qed

Note that the proof of Proposition \ref{prop:connected} uses the fact that $T$ is simple through  Lemma \ref{lem:(a,a)}(d) and hence requires  $f\geq 2$.

Putting together Lemma \ref {lem:gensubfield} and Proposition \ref{prop:connected}, we get the following corollary.

\begin{corollary}\label{cor:connected}
 Let $f\ge2$ and let $\GF(2^e)$ be the subfield generated by $\alpha$. 
\begin{itemize}
 \item[(a)] if $e$ is odd, or if $e$ is even and $\alpha^{2^{(e/2)}}\neq\alpha+1$, then $\Ga(f,\alpha)=m\Ga(e,\alpha)$, where $m=|T:\PSL(2,2^e)|^2$ and $\Ga(e,\alpha)$ is connected.
\item[(b)] if $e$ is even and $\alpha^{2^{(e/2)}}=\alpha+1$, then $\Ga(f,\alpha)$ has $|\PSL(2,2^e)|^{-2} |\PSL(2,2^f)|^3$ connected components, each isomorphic to $\Cos(\la \PSL(2,2^e),\nu\ra,H,H\ua\nu H)$, where $H=\PSL(2,2)$ and $\nu=\tau^{(e/2)}$.
\end{itemize}
\end{corollary}

We can now deal with the case $f=2$.
Take $\GF(4)=\{a+bi|a,b\in \GF(2), i^2=i+1\}$. , By Proposition \ref{prop:samegraph},
Construction \ref{cons} yields two graphs for $f=2$, namely $\Ga(2,0)$ and $\Ga(2,i)$.
\begin{corollary}\label{cor:f=2} The two graphs obtained by Construction \ref{cons} for $f=2$ are not connected. More precisely,
\begin{itemize} 
 \item[(a)]  $\Ga(2,0)\cong 200 \,K_{3,3}$, and 
\item[(b)] $\Ga(2,i)\cong 60 \, \mathcal{D}$ where $\mathcal{D}$ is the incidence graph of the Desargues configuration, called the Desargues graph (it is a double cover of the Petersen graph).
\end{itemize}
\end{corollary}
\proof
Consider first $\alpha=0$.
By Lemma \ref{lem:gensubfield},  $\Ga(2,0)\cong m\Ga(1,0)$, where $m=|\PSL(2,2^2):\PSL(2,2^1)|^2=100$. Part (a) follows from Proposition \ref{lem:f=1}.

Now assume $\alpha=i$. Then $\alpha^{2^{(f/2)}}=i^2=i+1=\alpha+1$, so part (b) of Proposition \ref{prop:connected} holds. Here $\ua=t_{1,i,0,1}$ and $\nu=\tau$. Thus $\Ga(2,i)$ has $|\PSL(2,2^2)|=60$ connected components, each containing $60/3=20$ vertices and isomorphic to $\Cos(\PGammaL(2,4),H,H\ua\tau H)$ where $H=\PSL(2,2)$.
There are only two arc-transitive cubic graphs on 20 vertices, the Desargues graph and the dodecagon (see \cite[p.148]{Biggs}). Since $\Ga(2,i)$ is bipartite by Proposition \ref {prop:generalities}, its connected components cannot be dodecagons, hence they are Desargues graphs. The Desargues graph has vertices the points and lines of the Desargues configuration, with two vertices adjacent if they form a flag (incident point-line pair) of the configuration. It is a double cover of the Petersen graph.
\qed

\section{Automorphism groups and isomorphisms for connected $\Ga(f,\alpha)$}

The remainder of this paper is concerned mainly with the connected graphs $\Ga(f,\alpha)$ given by Construction \ref{cons}, that is, we may assume from now on that $\alpha$ is a generator and, if $f$ is even, then  $\alpha^{2^{(f/2)}}\neq\alpha+1$ (see Corollary \ref{cor:connected}). By Lemma \ref{lem:f=1} and Corollary \ref{cor:f=2}, we may assume that $f\geq 3$.

In this section, we determine the full automorphism group $A$ of $\Ga=\Ga(f,\alpha)$ and the normaliser of $A$ in $\Sym(V\Ga)$. This will then enable us to determine a lower bound on the number of non-isomorphic such graphs, for a given $f$.

\begin{proposition}\label{prop:aut}
 Let $f\geqslant 3$ be an integer and $\alpha\in\GF(2^f)$.
Let $\Ga=\Ga(f,\alpha)$, $G$, $T$, $L$, $\pi$ be as in Construction \ref{cons} with $\Ga$ connected. 
The full automorphism group  of $\Ga$ is $A=G\times \la \sig\ra$, where $\sig$ is given by $(Lx)^\sig=L\pi x$ for all $x\in G$. In particular, $A$ does not depend on the choice of $\alpha$ and $\Ga$ is $(A,3)$-arc transitive but not $(A,4)$-arc-transitive.
Moreover, the stabiliser in $A$ of the vertex $L$ is $L\times \la \pi\sig\ra\cong D_{12}$.
\end{proposition}
\proof
Let $A$ be the full automorphism group of $\Ga$.
By Proposition \ref{prop:generalities}, $G\leqslant A$.
Define the map $\sig$ on $V\Ga$ by $(Lx)^\sig=L\pi x$ for all $x\in G$. This is a well defined bijection, since $\pi$ centralises $L$.
Consider an edge $\{Lg_1,Lg_2\}$, that is, $g_2g_1^{-1}\in Lg_\alpha L$. Its image under $\sig$ is  $\{L\pi g_1,L\pi g_2\}$. We have $\pi g_2(\pi g_1)^{-1}=\pi g_2g_1^{-1}\pi\in \pi Lg_\alpha L\pi=L\pi g_\alpha\pi L$. Recall that $g_\alpha=(u_\alpha,u_\alpha b)\pi$ and $u_\alpha b=b u_\alpha$, so 
$\pi g_\alpha\pi=(u_\alpha b,u_\alpha)\pi=(b,b)g_\alpha$. Thus $L\pi g_\alpha\pi L=Lg_\alpha L$, so $\{Lg_1,Lg_2\}^\sig$ is an edge, and $\sig\in A$. 
We now show that $\sig$ centralises $G$. Indeed, let $h\in G$ and $Lx\in V\Ga$, then $(Lx)^{h\sig}=(Lxh)^\sig=L\pi xh=(L\pi x)^h =(Lx)^{\sig h}$. Hence $\sig h=h\sig$, and $\sig\in \C_A(G)$.
Since $\Z(G)=1$, we have $\sig\notin G$. Also $\sig^2=1$. Therefore $R:=G\times \la \sig\ra \leqslant A$.
The stabiliser of $L\in V\Ga$ in $R$ is  $R_L=L\times \la \pi\sig\ra\cong S_3\times C_2\cong D_{12}$.

By Lemma \ref{lem:cosetgraph}(b), $\Ga$ is $(G,1)$-arc transitive, and so  is $(R,1)$-arc transitive. Tutte \cite{Tutte1,Tutte2} proved that the automorphism group of an arc-transitive finite graph with valency 3 acts regularly on $s$-arcs for some $s\leqslant 5$, and the stabiliser of a vertex has order $3.2^{s-1}$. Since $|R_L|=12$, $R$ acts regularly on the $3$-arcs of $\Ga$ (and hence is not transitive on 4-arcs).

Suppose $R<A$. Since  both $R$ and $A$ are transitive on $V\Ga$, the Orbit-Stabiliser Theorem implies that $R_L<A_L$, and so $A$ would act regularly on $s$-arcs for some $s=4$ or $5$. 
By Theorem 3 of \cite{DjoMil}, this is not possible. Hence $A=R$. 
\qed

\begin{definition}\label{def:isotau}{\rm 
Let $\Ga=\Ga(f,\alpha)$ (not necessarily connected).
We define $\bar\tau:V\Ga\rightarrow V\Ga:\,L(c,d)
\pi^\epsilon \mapsto \ L(c^\tau,d^\tau)\pi^\epsilon$ for each $(c,d)\in T^2, \epsilon\in\{0,1\}$, where  $\tau$ is as  defined in {\rm(\ref{tau})}.}
\end{definition}

\begin{lemma}\label{lem:isotau}
Let $\Ga=\Ga(f,\alpha)$  (not necessarily connected) and let  $\bar\tau$ be as in Definition \ref{def:isotau}.
Then $\bar\tau$ induces an isomorphism from $\Ga$ to $\Ga(f,\alpha^2)$.
Moreover $\la\bar\tau\ra\cong C_f$.
\end{lemma}
\proof
We have $\tau$, as  defined in {\rm(\ref{tau})}, in $\Aut(T)$. We denote by $\mu$ the element of $\Aut(G)$ defined by $(c,d)^\mu= (c^\tau,d^\tau)$ for all $(c,d)\in T^2$ and by $\pi^\mu=\pi$.
Then, since $\mu$ centralises $(a,a)$ and $(b,b)$, we have that $\mu\in\N_{\Aut\, G}(L)$. Thus we can use Lemma \ref{lem:cosetgraph}(f), with $\bar\mu:Lx\mapsto Lx^\mu$.
More precisely for  $(c,d)\in T^2, \epsilon\in\{0,1\}$, we have $(L(c,d)\pi^\epsilon)^{\bar\mu}=L(c,d)^\mu(\pi^\epsilon)^\mu=L(c^\tau,d^\tau)\pi^\epsilon$. 
Hence $\bar\mu=\bar\tau$ is a permutation of $V\Ga$ and induces an isomorphism from $\Ga=\Cos(G,H,H\ga H)$ to $\Cos(G,H,H\ga^\mu H)$ by Lemma \ref{lem:cosetgraph}(f).
Note that $\ga^\mu=((t_{1,\alpha,0,1}, t_{1,\alpha+1,0,1})\pi)^\mu$ (see (\ref{defg})), and so $\ga^\mu=((t_{1,\alpha,0,1})^\tau, (t_{1,\alpha+1,0,1})^\tau)\pi=(t_{1,\alpha^2,0,1}, t_{1,\alpha^2+1,0,1})\pi=g_{\alpha^2}$. Therefore $\Cos(G,H,H\ga^\mu H)=\Ga(f,\alpha^2)$. 

For $i\geq 1$, the permutation $\bar\tau^i$ of $V\Ga$ maps the coset $L(c,d)\pi^\epsilon$ onto $L(c^{\tau^i},d^{\tau^i})\pi^\epsilon$. Thus $\bar\tau$ has the same order as $\tau$, and so  $\la\bar\tau\ra\cong C_f$.
\qed

We now determine $\N_{\Sym(V\Gamma)}(A)$.
\begin{lemma}\label{lem:norm}
Let $\Ga=\Ga(f,\alpha)$ and $A$ be as in Proposition \ref{prop:aut}.
Then $\N_{\Sym(V\Gamma)}(A) = A \rtimes \la\bar\tau\ra\cong A . C_f$,
where $\bar\tau$ is as defined in Definition \ref{def:isotau}.
\end{lemma}

\proof Set $N:= \N_{\Sym(V\Gamma)}(A)$ and $N_0:= \la A, \bar\tau\ra$. 
We use the notation of Construction~\ref{cons}.
By Lemma \ref{lem:isotau}, $\bar\tau\in\Sym(V\Gamma)$. Moreover, it follows from the definitions
of $\bar\tau$ and $\sigma$ that 
${\bar\tau}^{-1}(c,d)\bar\tau= (c^\tau,d^\tau)$ for each $(c,d)\in T^2$, 
and $[\bar\tau,\sigma]=[\bar\tau,\pi]=1$. Thus 
$N_0 = A \rtimes \la\bar\tau\ra\leqslant N$ with $N_0/A\cong \la\bar\tau\ra\cong C_f$.

Since $T^2$ is a characteristic subgroup of $A$, each element of $N$
induces an automorphism of $T^2$, and we have a homomorphism 
$\varphi:N\rightarrow\Aut(T^2)$ with kernel
$K= \C_N(T^2)\leq  \C_{\Sym(V\Gamma)}(T^2)=\C$, say.
Now $K$ (and hence $\C$) contains $Z(A)=\la\sigma\ra\cong C_2$, 
which interchanges
the two orbits of $T^2$ in $V\Gamma$, and so the subgroup $\C^+$
of $\C$ stabilising each of the $T^2$-orbits setwise has index $2$ in $\C$.
The two $T^2$-orbits are the sets $\Delta_1$ and $\Delta_2$ of $L$-cosets
in $T^2$ and $T^2\ga$ respectively, and $L$ is the stabiliser in $T^2$
of the vertex $L$ of $\Delta_1$ and also the stabiliser in $T^2$ of 
the vertex $L\pi$ of $\Delta_2$. 
For $i=1,2$, let $S_i, L_i$ 
denote the permutation groups on $\Delta_i$ induced by $T^2$ and $L$ 
respectively. Then  by  Lemma \ref{lem:groupprop}(d), $\N_{S_i}(L_i)=L_i$ and by 
\cite[Theorem 4.2A(i)]{DM}, $\C_{\Sym(\Delta_i)}(S_i)\cong 
\N_{S_i}(L_i)/L_i = 1$.
Thus $\C^+=1$ and $K=\C=\la\sigma\ra$, of order 2.

Now $\varphi(N)$ contains the inner automorphism group $\varphi(T^2)$ of $T^2$,
and the quotient 
$\varphi(N)/\varphi(T^2)$ is contained in the outer automorphism group 
of $T^2$, which is isomorphic to
$\la\tau\ra\Wr\la\pi\ra$. Further, $\varphi(N)/\varphi(T^2)$ normalises 
$\varphi(A)/\varphi(T^2)$, which corresponds to the subgroup
$\la\pi\ra$ of $\la\tau\ra\Wr\la\pi\ra$, and so the subgroup of
$\la\tau\ra\Wr\la\pi\ra$ corresponding to $\varphi(N)/\varphi(T^2)$ 
lies in the normaliser of $\la\pi\ra$ in $\la\tau\ra\Wr\la\pi\ra$, namely 
$\la(\tau,\tau)\ra\times\la\pi\ra\cong 
C_f\times C_2$. On the other hand  $\varphi(N)/\varphi(T^2)$ contains
$\varphi(N_0)/\varphi(T^2)\cong \la\bar\tau\ra\times\la\pi\ra$. Thus 
equality holds and we conclude that $N=N_0$.
\qed

We are now able to determine a lower bound on the number of non-isomorphic connected graphs $\Ga(f,\alpha)$ for each $f$. 
They are obviously not isomorphic for different values of $f$, so in particular, it follows that there are infinitely many such graphs, as the lower bound is increasing with $f$. 

\begin{proposition}\label{numberisom}
 Let $f\geq 3$.
\begin{itemize}
\item[(a)] Let  $\Gamma(f,\alpha)$ and $\Gamma(f,\beta)$ be connected graphs.
Then $\Gamma(f,\alpha)\cong \Gamma(f,\beta)$ if and only if $\beta\in\{\alpha^{2^i}|0\leq i<f\}\cup\{\alpha^{2^i}+1|0\leq i<f\}$.
\item[(b)] The number of pairwise non-isomorphic connected graphs $\Ga$ obtained from Construction \ref{cons} is greater than $2^{f-2}/f$  if $f$ is odd and greater than $(2^{f-2}-2^{f/2-1})/f$ if $f$ is even.
\end{itemize}
\end{proposition}
\proof
Let  $\Ga=\Gamma(f,\alpha)$ and $\Gamma(f,\beta)$ be connected graphs produced by
Construction~\ref{cons}. 
By Corollary \ref{cor:connected}, $\alpha$ and $\beta$ are generators, and if $f$ is even then $\alpha^{2^{(f/2)}}\neq\alpha+1$ and  $\beta^{2^{(f/2)}}\neq\beta+1$.

Suppose  that $\psi$ is an isomorphism from $\Gamma(f,\alpha)$ to 
$\Gamma(f,\beta)$.  
Since $V\Gamma= V\Gamma(f,\beta)$, the isomorphism $\psi$ is an element of $\Sym(V\Ga)$ and since, by Proposition~\ref{prop:aut}, 
$\Aut(\Gamma(f,\alpha))=\Aut(\Gamma(f,\beta))=A$, it
follows that $\psi$ is an element of $\N_{\Sym(V\Gamma)}(A)$. 
By  Lemma~\ref{lem:norm},  $\N_{\Sym(V\Gamma)}(A)=A \rtimes \la\bar\tau\ra$.  Thus $\Gamma(f,\beta)$ is the image of $\Gamma(f,\alpha)$ under 
$\bar\tau^i$ for some $i$ such that $0\leq i<f$.
We have $\Gamma(f,\beta)=\Gamma(f,\alpha)^{\bar\tau^i}=\Ga(f,\alpha^{2^i})$ by Lemma \ref{lem:isotau}.
Therefore, by Proposition \ref{prop:samegraph}, $\beta = \alpha^{2^i}$ or $\alpha^{2^i}+1$, and so $\beta\in\{\alpha^{2^i}|0\leq i<f\}\cup\{\alpha^{2^i}+1|0\leq i<f\}$.

Suppose now that  $\beta\in\{\alpha^{2^i}|0\leq i<f\}\cup\{\alpha^{2^i}+1|0\leq i<f\}$. Then, by Proposition \ref{prop:samegraph}, $\Gamma(f,\beta)=\Gamma(f,\alpha^{2^i})$ for some $0\leq i<f$, which, by Lemma \ref{lem:isotau}, is equal to $\Gamma(f,\alpha)^{\bar\tau^i}$, where $\bar\tau^i$ is a graph isomorphism. Hence $\Gamma(f,\alpha)\cong \Gamma(f,\beta)$ and part (a) holds.

Let $\alpha$ be a generator such that, if $f$ is even,  $\alpha^{2^{(f/2)}}\neq\alpha+1$ .
We claim that the set $\{\alpha^{2^i}|0\leq i<f\}\cup\{\alpha^{2^i}+1|0\leq i<f\}$ has size $2f$. Notice first that all elements $x$ of this set are generators and do not satisfy the equation  $x^{2^{(f/2)}}\neq x+1$.
Suppose $\alpha^{2^i}=\alpha^{2^j}$ for some $i,j$ such that $0\leq i<j<f$, then  $\alpha^{2^i}=(\alpha^{2^i})^{2^{j-i}}$, contradicting  Lemma \ref{lem:alpha+1}(2) for the generator $\alpha^{2^i}$. Hence $\{\alpha^{2^i}|0\leq i<f\}$ and $\{\alpha^{2^i}+1|0\leq i<f\}$ both have size $f$. Now suppose $\alpha^{2^i}=\alpha^{2^j}+1$ for some $i,j$ such that $0\leq i<j<f$ (we can assume $i<j$ without loss of generality, because otherwise we just add $1$ to both sides of the equation).  Thus $\alpha^{2^i}=(\alpha^{2^i})^{2^{j-i}}+1$. Applying Lemma \ref{lem:alpha+1}(1) to the generator $\alpha^{2^i}$, we get that $f$ is even, $j-i=f/2$ and $\alpha^{2^i}=(\alpha^{2^i})^{2^{f/2}}+1$. However, since $\alpha^{2^i}$ does not satisfy the equation  $x^{2^{(f/2)}}\neq x+1$, this is a contradiction. Thus the claim is proved.

Suppose first $f$ is odd. Then $\Ga(f,\alpha)$ is connected if and only if $\alpha$ is a generator, by Corollary \ref{cor:connected}. By Lemma \ref{ff}, the number of generators of $\GF(2^{f})$ is strictly greater than $2^{f-1}$. By the claim and part (a), exactly $2f$ of those generators yield isomorphic graphs, thus the number of pairwise non-isomorphic connected graphs is greater than $2^{f-2}/f$.

Finally assume $f$ is even. Then $\Ga(f,\alpha)$ is connected if and only if $\alpha$ is a generator and  $\alpha^{2^{(f/2)}}\neq\alpha+1$, by Corollary \ref{cor:connected}.
By Lemma \ref{ffeven}, the number of such elements is greater than $2^{f/2}(2^{f/2-1}-1)$.  By the claim and part (a), exactly $2f$ of those generators yield isomorphic graphs, thus the number of pairwise non-isomorphic connected graphs is greater than $2^{f/2-1}(2^{f/2-1}-1)/f=(2^{f-2}-2^{f/2-1})/f$. 
\qed

We illustrate this result by considering the case $f=3$ where we obtain the first connected examples.
Take $\GF(8)=\{a+bj+cj^2|a,b,c\in \GF(2), j^3=j+1\}$.
For $f=3$ our construction yields four graphs with different edge-sets, namely $\Ga(3,0)$, $\Ga(3,j)$, $\Ga(3,j^2)$ and $\Ga(3,j^4)$, by Proposition \ref{prop:samegraph}.
\begin{corollary}\label{cor:f=3}
Up to isomorphism, Construction \ref{cons} for $f=3$ yields two graphs, one of which is connected. More precisely
\begin{itemize}
 \item[(a)]  $\Ga(3,0)\cong 14112\, K_{3,3}$, and
\item[(b)] $\Ga(3,j)\cong\Ga(3,j^2)\cong\Ga(3,j^4)$ is connected.
\end{itemize}
\end{corollary}
\proof
Consider first $\alpha=0$.
By Lemma \ref{lem:gensubfield},  $\Ga(3,0)\cong m\Ga(1,0)$, where $m=|\PSL(2,2^8):\PSL(2,2)|^2=84^2$. Part (a) now follows from Proposition \ref{lem:f=1}.
Now assume $\alpha=j$. By Proposition \ref{prop:connected}, $\Ga(3,j)$ is connected, and by Proposition \ref{numberisom}(a), $\Ga(3,j)\cong\Ga(3,j^2)\cong\Ga(3,j^4)$. 
\qed

For $f=4$ also, our construction yields just one connected graph and three disconnected ones, up to isomorphism.
Take $\GF(16)=\{a+bk+ck^2+dk^3|a,b,c,d\in \GF(2), k^4=k+1\}$.
\begin{corollary}\label{cor:f=4}
Up to isomorphism, Construction \ref{cons} for $f=4$ yields four graphs, one of which is connected. More precisely
\begin{itemize}
 \item[(a)]  $\Ga(4,0)=924800 \, K_{3,3}$, 
\item[(b)] for $\alpha\in\{k^5,k^{10}\}$, $\Ga(f,\alpha)\cong 277440 \, \mathcal{D}$, where $\mathcal{D}$ is the Desargues graph,
\item[(c)]  for $\alpha\in \{k,k^2,k^4,k^8\}$, $\Ga(4,\alpha)\cong \Ga(4,k)$ has $4080$ connected components, and
\item[(d)] for $\alpha$ a generator not in $\{k,k^2,k^4,k^8\}$,  $\Ga(4,\alpha)\cong\Ga(4,k^3)$ is connected.
\end{itemize}
\end{corollary}
\proof
Consider first $\alpha=0$.
By Lemma \ref{lem:gensubfield},  $\Ga(4,0)\cong m\Ga(1,0)$, where $m=|\PSL(2,16):\PSL(2,2)|^2=680^2$. Part (a) now follows from Proposition \ref{lem:f=1}.

The element $k^5$  generates $\GF(4)=\{0,1,k^5,k^{10}\}$, and so by Lemma \ref{lem:gensubfield}, $\Ga(4,k^5)\cong m\Ga(2,k^5)$, where $m=|\PSL(2,16):\PSL(2,4)|^2=68^2$. 
Now $\Ga(2,k^5)$ is $\Ga(2,i)$ from Corollary \ref{cor:f=2}, and so $\Ga(4,k^5)\cong 68^2.60 \, \mathcal{D}=277440\, \mathcal{D}$. Now $k^{10}=k^5+1$, and so by Proposition \ref{prop:samegraph}, $\Ga(4,k^5)=\Ga(4,k^{10})$. Thus part (b) holds.

Now assume $\alpha=k$. By Proposition \ref{numberisom}(a), $\Ga(4,k)\cong\Ga(4,k^2)\cong\Ga(4,k^4)\cong\Ga(4,k^8)$. Since $\alpha$ is a generator and  $\alpha^{2^{(f/2)}}=\alpha^4=\alpha+1$, by Proposition \ref{prop:connected}, $\Ga(4,k)$ has $|T|=4080$ connected components. Thus part (c) holds.

Finally assume $\alpha=k^3$. Then, by Proposition \ref{numberisom}(a),  $\Gamma(f,\beta)\cong \Gamma(f,k^3)$ if and only if $\beta\in\{\alpha^{2^i}|0\leq i<f\}\cup\{\alpha^{2^i}+1|0\leq i<f\}=\{k^3,k^6,k^{12},k^9\}\cup\{k^{14},k^{13},k^{11},k^7\}$, that is, if $\beta$ is any generator not in $\{k,k^2,k^4,k^8\}$.  Moreover, by Proposition \ref{prop:connected}, $\Ga(4,k^3)$ is connected since  $\alpha^4\neq \alpha+1$. Thus part (d) holds.
\qed

For $f=5$, the bound of Proposition \ref{numberisom} tells us that there are at least 2 non-isomorphic connected graphs obtained by Construction \ref{cons}. Actually there are $30$ generators, exactly $2f=10$ of them yielding isomorphic graphs, and so there are 3 pairwise non-isomorphic connected graphs for $f=5$.  

\section{Symmetry properties for connected $\Ga(f,\alpha)$}\label{sec:symm}
In this section, we  study the symmetry properties described in Tables \ref{def:trans} and \ref{def:localtrans} possessed by connected graphs $\Ga(f,\alpha)$. This includes a formal proof of Theorems \ref{thm:broadbrush} and \ref{thm:curiosity}. We start by defining the following five groups of automorphisms.

\begin{definition}\label{diffgroups}{\rm
We consider the following five subgroups of $A$, whose inclusions are given in Figure \ref{lattice}.
\begin{enumerate}
 \item $A=G\times\langle\sig\rangle$;
\item $A^+=T^2\rtimes \la \sig\pi\rangle$;
\item $G=T^2\rtimes \la\pi\ra$;
\item $M=T^2\times\la\sig\ra$;
\item $G^+=M^+=T^2$.
\end{enumerate}
}
\end{definition}

Note that $\sig\pi$ stabilises the bipartite halves of $\Ga(f,\alpha)$ setwise and  $T^2\rtimes \langle \sig\pi\rangle$ is maximal in $A$, hence it is $A^+$.
By Proposition \ref{prop:generalities}, $G^+=T^2$. Since $T^2$ is maximal in $M$ it follows that $M^+=T^2$.

\begin{figure}
\centering
$\begin{array}{rcccl}
           \; & \; & A     &\;       &\; \\
           \; & \diagup      & \mid & \diagdown &\; \\
           A^+     & \; & G & \;      &M \\
           \; & \diagdown & \mid & \diagup     &\; \\
           \; & \; & G^+   & \;  &\; \\
\end{array}$
\caption{Lattice} \label{lattice}
\end{figure}

We have the following results on $s$-arc transitivity.

\begin{proposition}\label{prop-arc}
Let $f\geq 3$, $\Ga(f,\alpha)$ be a connected graph as described in Construction \ref{cons}, and let $G,M,A,G^+,A^+$ be as in Definition \ref{diffgroups}. Then the following facts hold. 
\begin{enumerate}
\item $\Ga$ has girth at least 10.
 \item $\Ga$ is $(A,3)$-arc transitive but not $(A,4)$-arc transitive.
\item $\Ga$ is locally $(A^+,3)$-arc transitive but not locally $(A^+,4)$-arc transitive.
\item $\Ga$ is $(G,2)$-arc transitive but not $(G,3)$-arc transitive.
\item $\Ga$ is $(M,2)$-arc transitive but not $(M,3)$-arc transitive.
\item $\Ga$ is locally $(G^+,2)$-arc transitive but not locally $(G^+,3)$-arc transitive.
\end{enumerate}
\end{proposition}
\proof
See Proposition \ref{prop:aut} for (2).  Since $A^+_L=A_L$ has order $3.2^2$, we have that $\Ga$ is locally $(A^+,3)$-arc transitive but not locally $(A^+,4)$-arc transitive and (3) holds.

 By \cite[Theorem 2.1]{ConNed}, all the $3$-arc transitive finite graphs of girth up to 9 with valency 3 are known. The largest one has 570 vertices. By Theorem \ref{prop:generalities}, $|V\Ga|\geqslant 2^{6}(2^{6}-1)^2/3=84672$. Thus $\Ga$ has girth at least 10 and (1) holds.

Let $X\in\{G,G^+,M\}$. The stabiliser of the vertex ``$L$'' in $X$  is precisely $L$, acting as $S_3$ on its three neighbours. Therefore the stabiliser of a vertex is 2-transitive on its neighbours, and so $\Ga$ is locally $(X,2)$-arc transitive (see for instance Lemma 3.2 of \cite{GLP}). Since $G$ and $M$ are transitive on $V\Ga$, $\Ga$ is also $(G,2)$-arc transitive and $(M,2)$-arc transitive. Since $\girth(\Ga)>6$, the number of $3$-arcs  starting in $L$ is exactly $12$, and so $X_L$, which has order $6$, cannot be transitive on the 3-arcs starting in $L$.
Hence (4), (5) and (6) hold.
\qed

The lower bound of $10$ on the girth is an underestimate, but is sufficient for our purposes. For example, a computation using MAGMA \cite{magma} shows that, for $f=3$, the unique connected graph $\Ga(f,j)$ (see Corollary \ref{cor:f=3}) has girth 16 and for $f=4$, the girth of the unique connected graph $\Ga(3,k^3)$ (see Corollary \ref{cor:f=4})  is 30. 

\begin{question}
 Is the girth of the connected graphs obtained from Construction \ref{cons} unbounded?
\end{question}

Let $\Gamma$ be a graph of girth $g$. If $s\leqslant[\frac{g-1}{2}]$, then $\Ga$ is (locally) $(G,s)$-distance transitive if and only if $\Ga$ is (locally) $(G,s)$-arc transitive \cite[Lemma 7.2]{LDT}. Since $\Ga(f,\alpha)$  has girth at least 10 we have the following corollary to Proposition \ref{prop-arc}.

\begin{corollary}
Let $s\leq 4$, $\Ga=\Ga(f,\alpha)$ and $X\leqslant \Aut(\Ga)$. Then $\Ga$ is (locally) $(X,s)$-distance
transitive if and only if $\Ga$ is (locally) $(X,s)$-arc transitive
\end{corollary}

The following proposition determines, for each of the automorphism groups $X\in\{A,G,M\}$, whether $X$ is biquasiprimitive on vertices and whether $X^+$ is quasiprimitive on each bipartite half. Recall that $M^+=G^+$.
\begin{proposition}\label{prop-basic}
Let $f\geq 3$, $\Ga=\Ga(f,\alpha)$ be a connected graph described in Construction \ref{cons}, and let $G,M,A,G^+,A^+$ be as in Definition \ref{diffgroups}. Then $G$ is biquasiprimitive on $V\Ga$, while $M$ and $A$ are not  biquasiprimitive on $V\Ga$, and $A^+$ is quasiprimitive on each bipartite half, while $G^+$ is not.
\end{proposition}
\proof
We recall that $\sig$ centralises $G$.
Since $\pi$ (respectively $\sig\pi$) interchanges the two direct factors of $G^+$, $T^2$ is a minimal normal subgroup of $G$ and of $A^+$, and indeed is the unique minimal normal subgroup.
Since $T^2$ has two orbits on vertices, $G$ is biquasiprimitive on $V\Ga$.
Also $A^+$ is faithful and quasiprimitive on each of its orbits. 

Let $N=1\times T$, then $N$ is normal in $G^+$ and in $M$. Notice that $|N|=|T|=2^f(2^{2f}-1)$ is less than the number of vertices in each bipartite half. Hence $N$ is intransitive on each bipartite half and so $\Ga_N$ is nondegenerate. More precisely  $|V\Ga_N|=2^f(2^{2f}-1)/3$ with half the vertices in each bipartite half. Thus $G^+$ is not quasiprimitive on each bipartite half.

Now let $N'=\langle\sig\rangle$, then $N'$ is normal in $A$ and in $M$.  Obviously $N'$ (which has order $2$) is intransitive on each bipartite half and so $\Ga_{N'}$ is nondegenerate.  More precisely  $|V\Ga_{N'}|=|V\Ga|/2$. Thus $A$ and $M$ are not  biquasiprimitive on $V\Ga$.
 \qed

\begin{remark}{\rm
As mentioned in the introduction, if $G^+$ is not quasiprimitive on each bipartite half, which is the case here, then we can form a $G^+$-normal quotient and obtain a smaller locally $s$-arc-transitive graph. For $\Ga=\Ga(f,\alpha)$, we can quotient by $N=1\times T$. Now $G^+/N\cong T$, so this yields a locally $(T,2)$-arc transitive graph $\Ga_N$ such that $T$ has two orbits on vertices and the stabiliser of any vertex is isomorphic to $S_3$. Moreover, by  \cite[Theorem 1.1]{GLP},  $\Ga(f,\alpha)$ is a cover of this quotient.
Since $M$ normalises $N$, the group $M/N\cong T\times C_2$  also acts on $\Ga_N$. This action is vertex-transitive and hence $\Ga_N$ is $(M/N,2)$-arc transitive. In particular $\Ga_N$ is not semisymmetric. 
}\end{remark} 

In general, not all automorphisms of a quotient graph must arise from automorphisms of the original graph.

We can now prove our two main theorems.

\medskip \emph{Proof of Theorem \ref{thm:broadbrush}} 
By Proposition \ref{numberisom}(b), the number of non-isomorphic connected graphs $\Ga(f,\alpha)$ increases with $f$, and so there are an infinite number of such graphs. By Proposition \ref{prop-arc}(4) the graphs are $(G,2)$-arc transitive. Moreover, by Proposition \ref{prop-basic}, $G$ is biquasiprimitive on $V\Ga$ while  $G^+$ is not quasiprimitive on each bipartite half.
\qed

\emph{Proof of Theorem \ref{thm:curiosity}} By Proposition \ref{prop:aut}, $G$ has index 2 in $A=\Aut(\Ga)$, and by Proposition \ref{prop-arc}(2), $\Gamma$ is $(A,3)$-arc-transitive.  It follows from Proposition \ref{prop-basic} that $A$ is not biquasiprimitive on vertices and $A^+$ is quasiprimitive on each bipartite half.\qed

Next we verify that $G$ is indeed of the type given in \cite[Theorem 1.1(c)(i)]{bqp} as claimed in the introduction. First a definition:

\begin{definition}\label{def:ci}{\rm
A permutation group  $G\leqslant\Sym(\Omega)$ is biquasiprimitive of type (c)(i), as described in Theorem 1.1 of \cite{bqp}, if  $G$ is permutationally isomorphic to a group with the following properties.
\begin{itemize}
 \item[(a)] $|\Omega|=2m$ and the even subgroup $G^+\leqslant S_m\times S_m$ is equal to $\{(h,h^\varphi)|h \in H\}$, where $H\leqslant S_m$, $\varphi\in\Aut(H)$ and $\varphi^2$ is an inner automorphism of $H$.
 \item[(b)] $H$ has two intransitive minimal normal subgroups $R$ and $S$ such that $S=R^\varphi$, $R=S^\varphi$, and $R\times S$ is a transitive subgroup of $S_m$.
 \item[(c)] $\{(h,h^\varphi)|h \in R\times S\}$ is the unique minimal normal subgroup of $G$.
\end{itemize}
}\end{definition}

\begin{corollary}\label{cor:(c)(i)}
Let $f\geq 3$, $\Ga=\Ga(f,\alpha)$ be a connected graph described in Construction \ref{cons}, and $G$ also as in Construction \ref{cons}. Then $G\leqslant\Sym(V\Ga)$ is of type (c)(i), as described in Theorem 1.1 of \cite{bqp}.
\end{corollary}
\proof
By \cite[Theorem 1.2 and Proposition 4.1]{bqp}, a biquasiprimitive group acting 2-arc transitively on a bipartite graph must satisfy the conditions of (a)(i) or (c)(i) of Theorem 1.1 of \cite{bqp}. For groups satisfying (a)(i), the even subgroup is quasiprimitive on each bipartite half.
Since the permutation group induced by the action of $G^+=T^2$ on a bipartite half is not quasiprimitive, by Proposition \ref{prop-basic}, $G$ satisfies the conditions of (c)(i), and hence is of type (c)(i) as in Definition \ref{def:ci}.
More precisely, we have $m=|V\Ga|/2$, $H=T^2$, $\varphi=\pi$, $R=1\times T$, $S=T\times 1$, and $R\times S=T^2=G^+$.
\qed

The proof of Proposition \ref{prop-basic} shows that $\Ga$ is an $A$-normal double cover of its $A$-normal quotient $\Ga_{\la\sig\ra}$. 
We have $\{L,L\pi\}=L^{\la\sig\ra}$. A computation using MAGMA \cite{magma} shows that, when $f=3$, $L\pi$ is the unique vertex at maximal distance from $L$. In other words, $\Ga$ is antipodal with antipodal blocks of size $2$.

\begin{question}
 Let $f\geq 3$ and $\Ga=\Ga(f,\alpha)$ be a connected graph described in Construction \ref{cons}.
Is  $\Ga$ always antipodal with antipodal blocks of size $2$?
\end{question}

\end{document}